\documentclass[11pt]{article}
\usepackage{amsmath,amsfonts,amssymb,amsthm,graphicx,verbatim,subfig,mathtools,wasysym,tikz}
\pdfoutput=1
\usepackage[utf8]{inputenc}
\usepackage{anyfontsize}
\usepackage[all]{xy}
\usepackage[a4paper,left=25mm]{geometry}
\ifx\pdftexversion\undefined

\usepackage[a4paper,colorlinks,linkcolor=black,filecolor=black,citecolor=black,urlcolor=black,pdfstartview=FitH]{hyperref}
\fi

\font\sixbb=msbm6
\font\eightbb=msbm8
\font\twelvebb=msbm10 scaled 1095
\newfam\bbfam
\textfont\bbfam=\twelvebb \scriptfont\bbfam=\eightbb
                           \scriptscriptfont\bbfam=\sixbb
\def\bb{\fam\bbfam\twelvebb}
\newcommand{\Rea}{{\bb R}}

\newcommand{\Int}{{\bb Z}}

\newcommand{\cck}{\mathcal{K}}

\newtheorem{theorem}{\bf Theorem}[section]
\newtheorem{claim}[theorem]{\bf Claim}

\newtheorem{proposition}[theorem]{\bf Proposition}
\newtheorem{corollary}[theorem]{\bf Corollary}

\newtheorem{example}[theorem]{\bf Example}
\newcommand{\enp}{\begin{flushright} $\Box$ \end{flushright}}
\newcommand{\beq}[0]{\begin{equation}}
\newcommand{\enq}[0]{\end{equation}}

\newcommand{\pr}{{\rm Pr}}

\newcommand{\sym}{\mathbb{S}}

\newcommand{\vol}{{\rm vol}}

\title{Topological Connectivity of Random Permutation Complexes}
\author{Roy Meshulam\thanks{Department of Mathematics,
Technion, Haifa 32000, Israel. e-mail:
meshulam@technion.ac.il~. Supported by ISF grant 686/20.} \and Omer Moyal\thanks{Department of Mathematics,
Technion, Haifa 32000, Israel. e-mail: omer.moyal@campus.technion.ac.il~.}}

\begin{document}
\maketitle

\begin{abstract}
Let $\sym_n$ denote the symmetric group on $[n]=\{1,\ldots,n\}$ with the uniform probability measure.
For a permutation $\pi \in \sym_n$ let $X_{\pi}$ denote the simplicial complex
on the vertex set $[n]$ whose simplices are all $\{i_0,\ldots, i_m\} \subset [n]$
such that $i_0<\cdots<i_m$ and $\pi(i_0)<\cdots < \pi(i_m)$.
For $r \geq 0$ let $p_r(n)$ denote the probability that $X_{\pi}$ is not topologically $r$-connected
for $\pi \in \sym_n$. It is shown that for fixed $r \geq 0$ there exist constants $0<C_r, C_r' < \infty$ such that
\[
C_r \frac{(\log n)^r}{n} \leq  p_r(n) \leq C_r' \frac{(\log n)^{2r}}{n}.
\]
\end{abstract}

\section{Introduction}
\label{s:intro}


Let $\sym_n$ denote the symmetric group on $[n]=\{1,\ldots,n\}$ with the uniform probability measure.
For a permutation $\pi \in \sym_n$ let $\prec_{\pi}$ denote the partial order on $[n]$
given by $i \prec_{\pi} j$ if both $i <j$ and $\pi(i)<\pi(j)$.
Let $X_{\pi}$ denote the order complex of the poset $([n],\prec_{\pi})$, i.e. the simplicial complex
on the vertex set $[n]$ whose $m$-simplices are all $\{i_0,\ldots, i_m\} \subset [n]$
such that $i_0<\cdots<i_m$ and $\pi(i_0)<\cdots < \pi(i_m)$.
The following result had been proved in an equivalent formulation by Przytycki and Silvero \cite{PS18}, and independently but later by Chacholski, Levi and Meshulam \cite{CLM20}.
\begin{theorem}[\cite{PS18,CLM20}]
\label{t:hom}
For any $n \geq 1$ and $\pi \in \sym_n$ the complex $X_{\pi}$ is either contractible or is homotopy equivalent to a wedge of spheres
\begin{equation}
\label{e:hom}
X_{\pi}\simeq S^{k_1} \vee \cdots \vee S^{k_m}.
\end{equation}
Conversely, for any $k_1,\ldots,k_m \geq 0$ there exist an $n \geq 1$ and $\pi \in \sym_n$ such that (\ref{e:hom}) holds.
\end{theorem}
\begin{example}
Let
$\pi=\left(
\begin{array}{c}
1234567 \\
3254176
\end{array}
\right).
$
Then $X_{\pi}$ depicted in Figure 1 is homotopy equivalent to $S^1 \vee S^2$.
\end{example}

\begin{figure}
\begin{center}
\includegraphics[scale=0.3]{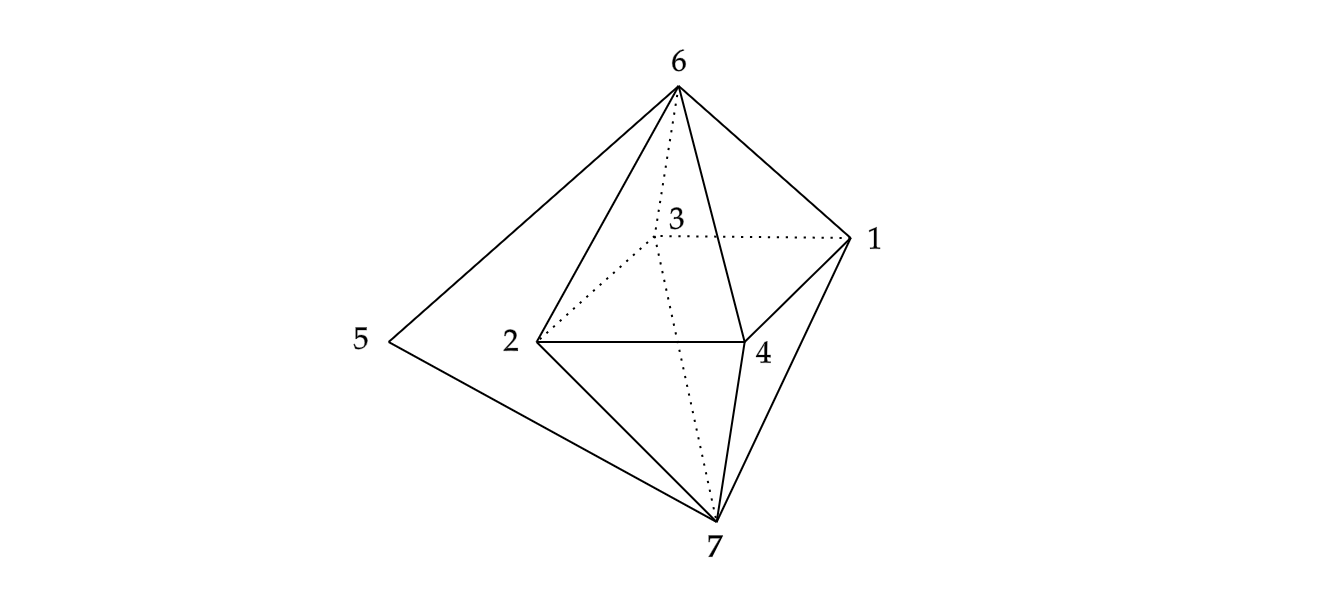}
\end{center}
\label{f:fig1}
\caption{$X_{\pi}$ for $\pi=\left(
\begin{array}{c}
1234567 \\
3254176
\end{array}
\right)
$}
\end{figure}

Let $r \geq 0$. Recall that a nonempty space $X$ is $r$-connected if for all $0 \leq i \leq r$,
any continuous map from the $i$-sphere $S^i$ to $X$ can be extended continuously to the $(i+1)$-th ball.
By convention, any nonempty space is $(-1)$-connected, while the empty set is not $(-1)$-connected.
Note that $X_{\pi}$ is $r$-connected if either $X_{\pi}$ is contractible, or if $r$ is strictly less than all $k_i$'s appearing
in the decomposition (\ref{e:hom}). In this paper we study the connectivity of $X_{\pi}$ for a random $\pi \in \sym_n$.
Let
\begin{equation*}
\label{e:def1}
p_r(n):=\pr\left[ \pi \in \sym_n: X_{\pi} \text{~is~not~} r\text{-connected} \right].
\end{equation*}
It is easy to check (see Claim \ref{c:pzero}) that $p_0(n)=\frac{2+o(1)}{n}$.
A much deeper result due to
Winkler \cite{Winkler85} asserts that asymptotically almost surely the $1$-skeleton of $X_{\pi}$ has diameter $3$.
Let $c_0=0$ and for $r \geq 1$ let $c_r=\sum_{i=0}^{r-1} \frac{1}{i!}$.
Our main result is the following
\begin{theorem}
\label{t:maint}
Let $r \geq 0$ and $n \geq 1$. Then
\begin{itemize}
\item[(i)]
\begin{equation*}
\label{e:mainu}
 p_r(n) \leq 40^{r+1} \frac{(\log 3n)^{2r}}{n}.
\end{equation*}
\item[(ii)]
\begin{equation*}
\label{e:mainl}
p_r(n) \geq \frac{1}{r!}\frac{(\log n)^r}{n}-c_r\frac{(\log n)^{r-1}}{n}.
\end{equation*}
\end{itemize}
\end{theorem}
For the proof of the upper bound it will be convenient to use an equivalent model for random permutation complexes (see e.g. \cite{Winkler85}). Let $\leq$ be the standard partial order on $\Rea^2$, i.e. $q=(a,b) \leq q'=(a',b')$
if $a \leq a'$ and $b \leq b'$. For $n \geq 1$ and an $n$-tuple $q=(q_1,\ldots,q_n) \in (\Rea^2)^n$ let
$Y_q$ denote the order complex of $\{q_1,\ldots,q_n\}$, i.e. the simplicial complex on the vertex set $[n]$ whose $m$-simplices are the subsets $\{i_0,\ldots,i_m\}$ such that
$q_{i_0}  <  \cdots <  q_{i_m}$. Let $\Omega_n=([0,1]^2)^n$ denote the probability space of $n$-tuples $q=(q_1,\ldots,q_n)$
where the $q_i$'s are picked independently and uniformly from $[0,1]^2$. Let $\Omega_n'$ denote the subspace of all
$q=(q_1,\ldots,q_n) \in \Omega_n$ with $q_i=(a_i,b_i)$ such that $|\{a_1,\ldots,a_n\}|=|\{b_1,\ldots,b_n\}|=n$.
For $q \in \Omega_n'$ as above let $\alpha,\beta \in \sym_n$ be the unique permutations such that
$a_{\alpha(1)}<\cdots<a_{\alpha(n)}$ and $b_{\beta(1)}<\cdots<b_{\beta(n)}$, and let
$\varphi(q)=\beta^{-1}\alpha \in \sym_n$.
 Extending $\varphi$ arbitrarily to the whole of $\Omega_n$, it is clear that
$\pr_{\Omega_n}\left[\varphi^{-1}(\pi)\right]=\frac{1}{n!}=\pr_{\sym_n}[\left\{\pi\right\}]$.
As $\pr_{\Omega_n}[\Omega_n']=1$ and $Y_q \cong X_{\varphi(q)}$ for all $q \in \Omega_n'$, it follows that
\begin{equation*}
\label{e:def2}
p_r(n)= \pr\left[q \in \Omega_n: Y_q \text{~is~not~} r\text{-connected}\right].
\end{equation*}

The paper is organized as follows. In Section \ref{s:nerve} we give a simple nerve type sufficient condition (Proposition \ref{p:suffc}) for $r$-connectivity of
$X_{\pi}$. In Section \ref{s:integrals} we establish an upper bound on a certain integral related to $p_r(n)$
(Proposition \ref{p:mainest1}). The results of Sections \ref{s:nerve} and \ref{s:integrals} are the key ingredients in the proof
of Theorem \ref{t:maint}(i) given in Section \ref{s:rconn}. In Section \ref{s:lbnd} we prove Theorem \ref{t:maint}(ii)
using a homotopy decomposition of $X_{\pi}$ (Proposition \ref{p:homdec}).
We conclude in Section \ref{s:con} with some remarks and open problems.

\section{A Sufficient Condition for $r$-Connectivity}
\label{s:nerve}
Let $K$ be a finite simplicial complex and let $\cck=\{K_j\}_{j=1}^m$ be a family of subcomplexes
of $K$ such that $\bigcup_{j=1}^m K_j=K$. For a nonempty $J \subset [m]$ let $K_{J}=\bigcap_{j \in J} K_j$.
The {\em nerve} $N(\cck)$ of $\cck$ is the simplicial complex on the vertex set $[m]$
whose simplices are all $J \subset [m]$ such that $K_{J} \neq \emptyset$.
As usual, let $\pi_0(K)$ denote the set of connected components of $K$ and for $i \geq 1$ let $\pi_i(K)$ denote the $i$-th homotopy group of $K$.
We shall need Bj\"{o}rner's version of the nerve theorem.
\begin{theorem}[\cite{Bjorner03}]
\label{t:nerve}
Let $r \geq 0$. Assume that $K_{J}$ is either empty or $(r-|J|+1)$-connected for all $\emptyset \neq J \subset [m]$ such that
$|J| \leq r+1$. Then
$\pi_i(K) \cong \pi_i(N(\cck))$ for all $0 \leq i \leq r$.
In particular, $K$ is $r$-connected iff $N(\cck)$ is $r$-connected.
\end{theorem}
\noindent
Fix $q=(q_1,\ldots,q_n) \in \Omega_n'$ where $q_i=(a_i,b_i)$ and let
\[
M(q)=\left\{i \in [n]: q_i \text{~is~a~minimal~element~of~}\{q_1,\ldots,q_n\}\right\}.
\]
For a subset $S \subset [n]$ let
$Y_q[S]=\{\sigma \in Y_q: \sigma \subset S\}$ denote the induced subcomplex of $Y_q$ on the vertex set $S$.
For $1 \leq i \neq  j \leq n$ let $q_i \vee q_j= (\max\{a_i,a_j\},\max\{b_i,b_j\})$ and
let
\[
Y_{q,ij}=Y_q\left[\left\{\ell: q_{\ell} \geq q_i \vee q_j\right\}\right].
\]
\begin{proposition}
\label{p:suffc}
Let $r\geq 1$. If $Y_{q,ij}$ is $(r-1)$-connected for all distinct $i,j \in M(q)$ then
$Y_q$ is $r$-connected.
\end{proposition}
\noindent
{\bf Proof.} Let $K=Y_q$ and for $i \in M(q)$ let
$K_i=Y_q\left[\{j : q_j \geq q_i\}\right]$.
Then $K=\bigcup_{i \in M(q)} K_i$.
Let $\emptyset \neq J \subset M(q)$ such that $|J| \leq r+1$. If $J=\{j\}$ is a singleton, then
$K_J=K_j$ is a cone with apex $j$, hence contractible and therefore $\ell$-connected for all $\ell$.
Suppose on the other hand that $|J| \geq 2$. Let $i,j \in J$ satisfy
$a_i=\max \{a_t: t \in M(q)\}$ and $b_j=\max \{b_t: t \in M(q)\}$. Then
$K_J=Y_{q,ij}$. By assumption $K_J$ is $(r-1)$-connected
and thus $(r-|J|+1)$-connected. Therefore $\{K_i\}_{i \in M(q)}$ satisfies the conditions of Theorem \ref{t:nerve}.
As $K_J \neq \emptyset$ for all $J \subset M(q)$, it follows that the nerve $N(\{K_i\}_{i\in M(q)})$ is a simplex and hence $r$-connected. Therefore
$K=Y_q$ is $r$-connected as well.
{\enp}

\section{A Key Estimate}
\label{s:integrals}
For $k \geq 1$ let $Q_k=\{(x_1,\ldots,x_{k}) \in \Rea^{k}: x_i \geq 0, \sum_{i=1}^{k} x_i \leq 1\}$
denote the $k$-simplex in $\Rea^k$ and let $\Delta_k=\{(x_1,\ldots,x_{k+1}) \in \Rea^{k+1}: x_i \geq 0, \sum_{i=1}^{k+1} x_i=1\}$
denote the $k$-simplex in $\Rea^{k+1}$. Recall (see e.g. Exercise 13, Chapter 10 in \cite{Rudin76})  that for any integers $\alpha_1,\ldots,\alpha_{k+1} \geq 0$
\begin{equation}
\label{e:intsim}
\begin{split}
&\int_{x \in \Delta_k} x_1^{\alpha_1}\cdots x_{k+1}^{\alpha_{k+1}} dx=\sqrt{k+1} \,
\int_{x \in Q_k} x_1^{\alpha_1}\cdots x_{k}^{\alpha_{k}} (1-\sum_{i=1}^k x_i)^{\alpha_{k+1}} dx \\
&=\frac{\sqrt{k+1}\, \alpha_1! \cdots \alpha_{k+1}!}{(k+\alpha_1+\cdots+\alpha_{k+1})!}.
\end{split}
\end{equation}
In particular, the $k$-dimensional volume of $\Delta_k$
is $\vol_k(\Delta_k)=\sqrt{k+1}\, \vol_k(Q_k)=\frac{\sqrt{k+1}}{k!}$.
 For $p \geq 1$ let $A_{k,p}=\{(\alpha_1,\ldots,\alpha_p) \in \Int^p: \alpha_i \geq 0, \sum_{i=1}^p \alpha_i=k\}$. For $k,\ell \geq 0$ define
\[
I(k,\ell)=\vol_2(\Delta_2)^{-2} \int_{(x,y) \in \Delta_2^2} \left(x_1y_3+x_2(y_2+y_3)+x_3(y_1+y_2)\right)^k (x_3y_3)^{\ell} dx \, dy.
\]
The main result of this section is the following estimate.
\begin{proposition}
\label{p:mainest1}
\begin{equation}
\label{e:mainest1}
I(k,\ell) \leq
\frac{40}{k+2} \binom{k+\ell+2}{2}^{-1}\binom{k+\ell}{k}^{-1}(1+\log (k+1)).
\end{equation}
\end{proposition}
\noindent
The proof of Proposition \ref{p:mainest1} depends on some preliminary observations.
For integers $N,t \geq 0$ such that $N \geq 2t$ let
\[
F(N,t)=\binom{N}{\lceil\frac{N}{2}\rceil+t} \sum_{i=0}^t \binom{N}{\lceil\frac{N}{2}\rceil+i}^{-1}.
\]
\begin{claim}
\label{c:fnt1}
\begin{equation*}
\label{e:fnt}
F(N,t)\leq 1+\min\left\{t,\frac{N}{t}\right\}.
\end{equation*}
\end{claim}
\noindent
{\bf Proof.} The bound $F(N,t) \leq 1+t$ is clear.
For the other inequality observe that if $0 \leq i \leq t$ then
\begin{equation}
\label{e:nti3}
\begin{split}
&\binom{N}{\lceil\frac{N}{2}\rceil+t}\binom{N}{\lceil\frac{N}{2}\rceil+i}^{-1} =
\frac{\left(\lceil\frac{N}{2}\rceil+i\right)! \left(\lfloor\frac{N}{2}\rfloor-i\right)!}{\left(\lceil\frac{N}{2}\rceil+t\right)!
\left(\lfloor\frac{N}{2}\rfloor-t\right)!} \\
&=\prod_{j=0}^{t-i-1}\frac{\lfloor\frac{N}{2}\rfloor-i-j}{\lceil\frac{N}{2}\rceil+t-j}
\leq \left(\frac{\frac{N}{2}-i}{\frac{N}{2}+t}\right)^{t-i}=\left(1-\frac{t+i}{\frac{N}{2}+t}\right)^{t-i} \\
&\leq \exp \left(-\frac{(t+i)(t-i)}{\frac{N}{2}+t}\right) \leq \exp\left(-\frac{t^2}{N}\right) \cdot
\exp\left(\frac{i^2}{N}\right).
\end{split}
\end{equation}
\noindent
Moreover,
\begin{equation}
\label{e:nti4}
\begin{split}
&\sum_{i=0}^{t-1}\exp\left(\frac{i^2}{N}\right) \leq \int_{0}^t \exp\left(\frac{x^2}{N}\right) dx
=\int_{0}^t \left(\sum_{j=0}^{\infty} \frac{x^{2j}}{N^j j!}\right) dx \\
&= \sum_{j=0}^{\infty} \frac{1}{N^j j!} \int_{0}^t x^{2j} dx
=\sum_{j=0}^{\infty} \frac{t^{2j+1}}{N^j j! (2j+1)} \\
&\leq \frac{N}{t} \sum_{j=0}^{\infty} \frac{t^{2(j+1)}}{N^{j+1} (j+1)!} =\frac{N}{t}
\left(\exp\left(\frac{t^2}{N}\right) -1\right).
\end{split}
\end{equation}
\noindent
Combining (\ref{e:nti3}) and (\ref{e:nti4}) it follows that
\begin{equation*}
\begin{split}
F(N,t)&=1+\binom{N}{\lceil\frac{N}{2}\rceil+t} \sum_{i=0}^{t-1} \binom{N}{\lceil\frac{N}{2}\rceil+i}^{-1}  \\
&\leq 1+ \exp\left(-\frac{t^2}{N}\right) \sum_{i=0}^{t-1}\exp\left(\frac{i^2}{N}\right) \\
&\leq 1+\exp\left(-\frac{t^2}{N}\right) \cdot \frac{N}{t}
\left(\exp\left(\frac{t^2}{N}\right) -1\right)<1+\frac{N}{t}.
\end{split}
\end{equation*}
{\enp}
\begin{corollary}
\label{c:fnt2}
For any $\ell, j \geq 0$
\begin{equation}
\label{e:fnt2}
\binom{2\ell+j}{\ell+j} \sum_{i=0}^j \binom{2\ell+j}{\ell+i}^{-1} \leq \frac{10 (\ell+j+1)}{j+1}.
\end{equation}
\end{corollary}
\noindent
{\bf Proof.} Write $N=2\ell+j$. Then
\begin{equation*}
\label{e:fellj}
\begin{split}
&\sum_{i=0}^j \binom{2\ell+j}{\ell+i}^{-1} = \sum_{i=0}^j \binom{N}{\lceil\frac{N}{2}\rceil -\lceil\frac{j}{2}\rceil+i}^{-1}\\
&= \sum_{i=0}^{\lceil\frac{j}{2}\rceil-1} \binom{N}{\lceil\frac{N}{2}\rceil -\lceil\frac{j}{2}\rceil+i}^{-1}
+\sum_{i=\lceil\frac{j}{2}\rceil}^j \binom{N}{\lceil\frac{N}{2}\rceil -\lceil\frac{j}{2}\rceil+i}^{-1} \\
&=\sum_{i=1}^{\lceil\frac{j}{2}\rceil} \binom{N}{\lfloor \frac{N}{2}\rfloor+i}^{-1}
+\sum_{i=0}^{\lfloor\frac{j}{2}\rfloor} \binom{N}{\lceil \frac{N}{2}\rceil+i}^{-1} \\
&=2\sum_{i=0}^{\lfloor\frac{j}{2}\rfloor} \binom{N}{\lceil \frac{N}{2}\rceil+i}^{-1} -\frac{1+(-1)^N}{2} \binom{N}{\lceil \frac{N}{2} \rceil}^{-1}.
\end{split}
\end{equation*}
Hence
\begin{equation*}
\begin{split}
&\binom{2\ell+j}{\ell+j} \sum_{i=0}^j \binom{2\ell+j}{\ell+i}^{-1}
=\binom{N}{\lceil\frac{N}{2}\rceil+\lfloor \frac{j}{2} \rfloor}
\sum_{i=0}^j \binom{2\ell+j}{\ell+i}^{-1} \\
&\leq 2\binom{N}{\lceil\frac{N}{2}\rceil+\lfloor \frac{j}{2} \rfloor}\sum_{i=0}^{\lfloor\frac{j}{2}\rfloor} \binom{N}{\lceil \frac{N}{2}\rceil+i}^{-1}
=2F(N,\lfloor\frac{j}{2}\rfloor) \\
&\leq 2+2\min \left\{\lfloor\frac{j}{2}\rfloor, \frac{N}{\lfloor\frac{j}{2}\rfloor}\right\} =
2+2\min \left\{\lfloor\frac{j}{2}\rfloor, \frac{2\ell+j}{\lfloor\frac{j}{2}\rfloor}\right\} \\
&\leq \frac{10(\ell+j+1)}{j+1}.
\end{split}
\end{equation*}
{\enp}

\begin{claim}
\label{c:intc1}
For integers $\beta_1,\beta_2,\beta_3 \geq 0$
\begin{equation}
\label{e:intc1}
\vol_2(\Delta_2)^{-1}\int_{y \in \Delta_2} y_3^{\beta_1}(y_2+y_3)^{\beta_2} (y_1+y_2)^{\beta_3} dy \leq
\frac{2\beta_1!\beta_3!}
{(\beta_2+1)(\beta_1+\beta_3+1)!}.
\end{equation}
\end{claim}
\noindent
{\bf Proof.}
Applying a change of variables $\phi:\{(z_1,z_2):0 \leq z_1\leq z_2 \leq 1\} \rightarrow \Delta_2$
 given by $\phi(z_1,z_2)=(1-z_2,z_2-z_1,z_1)$ with Jacobian $|J_{\phi}(z_1,z_2)|=\sqrt{3}$, we obtain
\begin{equation*}
\label{e:inty}
\begin{split}
&\vol_2(\Delta_2)^{-1}\int_{y \in \Delta_2} y_3^{\beta_1} (y_2+y_3)^{\beta_2} (y_1+y_2)^{\beta_3} dy \\
&=\sqrt{3}\,\vol_2(\Delta_2)^{-1}
\int_{0 \leq z_1 \leq z_2 \leq 1} z_1^{\beta_1}z_2^{\beta_2}(1-z_1)^{\beta_3} dz_1 \, dz_2 \\
&=2\int_{z_1=0}^1 z_1^{\beta_1} (1-z_1)^{\beta_3} \left(\int_{z_2=z_1}^1 z_2^{\beta_2} dz_2\right)
dz_1 \\
& \leq \frac{2}{\beta_2+1}\int_{z_1=0}^1 z_1^{\beta_1} (1-z_1)^{\beta_3}
dz_1 \\
&=\frac{2\beta_1!\beta_3!}
{(\beta_2+1)(\beta_1+\beta_3+1)!}.
\end{split}
\end{equation*}
{\enp}
\noindent
{\bf Proof of Proposition \ref{p:mainest1}.}
\begin{equation*}
\label{e:morec1}
\begin{split}
&I(k,\ell)= \vol_2(\Delta_2)^{-2}
\int_{(x,y) \in \Delta_2^2} \left(x_1y_3+x_2(y_2+y_3)+x_3(y_1+y_2)\right)^k (x_3y_3)^{\ell} dx \, dy \\
&= \sum_{\alpha \in A_{k,3}}
\binom{k}{\alpha_1;\alpha_2;\alpha_3} \vol_2(\Delta_2)^{-2}\int_{(x,y) \in \Delta_2^2}  (x_1y_3)^{\alpha_1} (x_2(y_2+y_3))^{\alpha_2}
(x_3(y_1+y_2))^{\alpha_3} (x_3y_3)^{\ell} dx \, dy \\
&=\sum_{\alpha \in A_{k,3}}
\frac{k!}{\alpha_1! \alpha_2! \alpha_3!} \left(\vol_2(\Delta_2)^{-1}\int_{x \in \Delta_2}
x_1^{\alpha_1}x_2^{\alpha_2}x_3^{\alpha_3+\ell} dx\right) \left(\vol_2(\Delta_2)^{-1} \int_{y \in \Delta_2}
y_3^{\alpha_1+\ell}(y_2+y_3)^{\alpha_2}
(y_1+y_2)^{\alpha_3} dy \right)\\
&\stackrel{(a)}{\leq} \sum_{\alpha \in A_{k,3}}
\frac{k!}{\alpha_1! \alpha_2! \alpha_3!}  \cdot
\left(\frac{2\alpha_1! \alpha_2! (\alpha_3+\ell)!}{(k+\ell+2)!} \right)
\cdot
\left(\frac{2(\alpha_1+\ell)! \alpha_3!}{(\alpha_2+1) (\alpha_1+\alpha_3+\ell+1)!}\right) \\
&=\frac{4k!}{(k+\ell+2)!} \sum_{\alpha \in A_{k,3}}
\frac{(\alpha_1+\ell)!(\alpha_3+\ell)!}{(\alpha_2+1)(\alpha_1+\alpha_3+\ell+1)!} \\
&=\frac{4k!}{(k+\ell+2)!} \sum_{j=0}^k
\frac{1}{(k-j+1)(j+\ell+1)!}
\sum_{(\alpha_1,\alpha_3) \in A_{j,2}}
(\alpha_1+\ell)!(\alpha_3+\ell)! \\
&=\frac{4k!}{(k+\ell+2)!} \sum_{j=0}^k
\frac{(2\ell+j)!}{(k-j+1)(j+\ell+1)!}
\sum_{i=0}^j \binom{2\ell+j}{\ell+i}^{-1} \\
&=\frac{4k!}{(k+\ell+2)!} \sum_{j=0}^k
\frac{\ell! (\ell+j)!}{(k-j+1)(j+\ell+1)!}
\binom{2\ell+j}{\ell+j}\sum_{i=0}^j \binom{2\ell+j}{\ell+i}^{-1} \\
&\stackrel{(b)}{\leq}\frac{4k!\ell!}{(k+\ell+2)!} \sum_{j=0}^k
\frac{1}{(k-j+1)(j+\ell+1)} \cdot \frac{10 (\ell+j+1)}{j+1} \\
&=\frac{40 k!\ell!}{(k+2)(k+\ell+2)!} \sum_{j=0}^k
\left(\frac{1}{k-j+1}+ \frac{1}{j+1}\right) \\
&= \frac{40}{k+2} \binom{k+\ell+2}{2}^{-1} \binom{k+\ell}{k}^{-1} \sum_{j=1}^{k+1}
\frac{1}{j} \\
 &\leq \frac{40}{k+2} \binom{k+\ell+2}{2}^{-1}\binom{k+\ell}{k}^{-1}(1+\log (k+1)),
\end{split}
\end{equation*}
where (a) follows from (\ref{e:intsim}) and (\ref{e:intc1}), and (b) follows from (\ref{e:fnt2}).
{\enp}

\section{The Upper Bound}
\label{s:rconn}

In this section we prove Theorem \ref{t:maint}(i). We first obtain a recursive upper bound on $p_r(n)$.
\begin{proposition}
\label{p:recursion}
Let $r \geq 1$ and $n \geq 3$. Then
\begin{equation}
\label{e:recursion}
p_r(n) \leq 20 \left(1+\log n \right) \left( \frac{1}{n}+\sum_{k=0}^{n-3} \frac{p_{r-1}(n-2-k)}{k+2} \right). \\
\end{equation}
\end{proposition}
\noindent
{\bf Proof.} For $1 \leq i<j \leq n$ and a subset $D \subset [n] \setminus \{i,j\}$ let
 $R_{i,j,D}$ denote the set of all $q=(q_1,\ldots,q_n) \in \Omega_n$
such that
\begin{itemize}
\item[(C1)]
$i,j \in M(q)$.
\item[(C2)]
$\{t \in [n] \setminus \{i,j\}: q_t \ngeq q_i \vee q_j\}=D$.
\end{itemize}
Fix $c=(a_1,a_2,b_1,b_2) \in [0,1]^4$ such that $a_1<a_2$ and $b_2<b_1$,
and let
\begin{equation*}
\begin{split}
E_-(c)&=
[0,a_1]\times [0,b_1] \cup [0,a_2] \times [0,b_2] \setminus \{(a_1,b_1),(a_2,b_2)\}, \\
E_+(c)&=[a_2,1] \times [b_1,1], \\
E(c)&=[0,1]^2 \setminus \left(E_-(c) \cup E_+(c)\right).
\end{split}
\end{equation*}

\begin{figure}
\begin{center}
\includegraphics[scale=0.34]{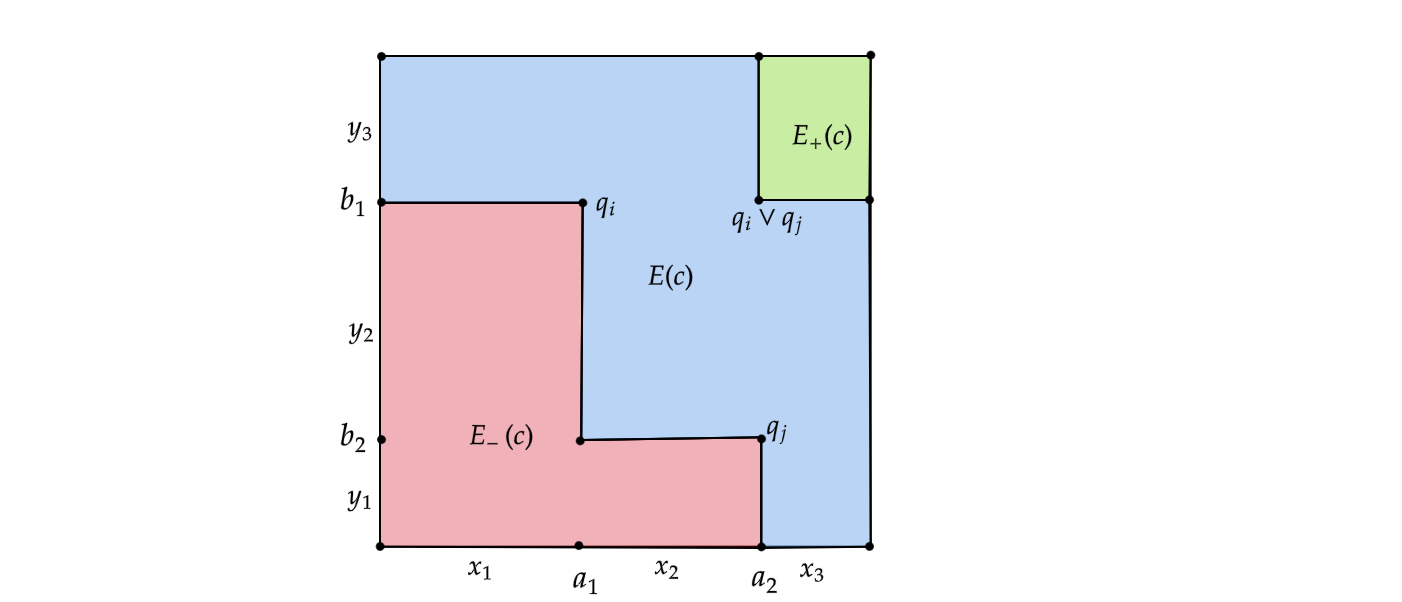}
\end{center}
\caption{$E_{-}(c), E(c), E_+(c)$}
\end{figure}
\noindent
Suppose now that $q=(q_1,\ldots,q_n) \in \Omega_n$ satisfies $q_i=(a_1,b_1), q_j=(a_2,b_2)$, see Figure 2.
Then $q \in R_{i,j,D}$ iff
 $q_t \in E(c)$ for $t \in D$,
and
$q_t \in E_+(c)$ for $t \in [n] \setminus (\{i,j\} \cup D)$.
Writing
\[
\begin{split}
x&=(x_1,x_2,x_3)=(a_1,a_2-a_1,1-a_2), \\
y&=(y_1,y_2,y_3)=(b_2,b_1-b_2,1-b_1),
\end{split}
\]
 it follows that
\begin{equation}
\label{qiin3}
\begin{split}
&\pr\left[~q=(q_1,\ldots,q_n) \in R_{i,j,D} ~|~q_i=(a_1,b_1), q_j=(a_2,b_2) \right] \\
&=\vol_2\left(E(c)\right)^{|D|} \vol_2\left(E_+(c)\right)^{n-2-|D|} \\
&=\left(a_1(1-b_1)+(a_2-a_1)(1-b_2)+(1-a_2)b_1\right)^{|D|} \left((1-a_2)(1-b_1)\right)^{n-2-|D|} \\
&= \left(x_1y_3+x_2(y_2+y_3)+x_3(y_1+y_2)\right)^{|D|} \left(x_3y_3\right)^{n-2-|D|}.
\end{split}
\end{equation}
\noindent
Combining (\ref{qiin3}) and (\ref{e:mainest1}) and noting that the probability that $q_i$ and $q_j$ are incomparable is $\frac{1}{2}$
we obtain
\begin{equation*}
\label{e:prnk}
\begin{split}
&\pr\left[R_{i,j,D}\right]=\frac{1}{2} \vol_2(\Delta_2)^{-2}
\int_{(x,y) \in \Delta_2^2} \left(x_1y_3+x_2(y_2+y_3)+x_3(y_1+y_2)\right)^{|D|} \left(x_3y_3\right)^{n-2-|D|} dx \, dy \\
&=\frac{1}{2} I(|D|,n-2-|D|) \leq \frac{20}{|D|+2} \binom{n}{2}^{-1}\binom{n-2}{|D|}^{-1}\left(1+\log (|D|+1)\right).
\end{split}
\end{equation*}
Proposition \ref{p:suffc} implies that if $Y_q$ is not $r$-connected then there exist $1 \leq i < j \leq n$ and
$D \subset [n]\setminus \{i,j\}$ such that
$q \in R_{i,j,D}$ and the induced complex $Y_q\left[[n]\setminus (D \cup\{i,j\})\right]$ is not $(r-1)$-connected.
Therefore by the union bound
\begin{equation*}
\label{e:prprm}
\begin{split}
&p_r(n) \leq
\sum_{1 \leq i < j \leq n} \sum_{D \subset [n] \setminus \{i,j\}}
\pr\left[q \in R_{i,j,D}:Y_q\left[[n]\setminus (D \cup\{i,j\})\right] \text{~is~not~} (r-1)\text{-connected}\right] \\
&=
\sum_{1 \leq i < j \leq n} \sum_{D \subset [n] \setminus \{i,j\}}
\pr\left[R_{i,j,D}\right] p_{r-1}(n-2-|D|) \\
&\leq \binom{n}{2} \sum_{k=0}^{n-2} \binom{n-2}{k} \left( \frac{20}{k+2} \binom{n}{2}^{-1}\binom{n-2}{k}^{-1}\left(1+\log n \right) \right)
\cdot  p_{r-1}(n-2-k) \\
&=20 \left(1+\log n \right) \sum_{k=0}^{n-2} \frac{p_{r-1}(n-2-k)}{k+2} \\
&=20 \left(1+\log n \right) \left( \frac{1}{n}+\sum_{k=0}^{n-3} \frac{p_{r-1}(n-2-k)}{k+2} \right).
\end{split}
\end{equation*}
\enp
We also need the following observation.
\begin{claim}
\label{c:pzero}
(i) $\frac{1}{n} \leq p_0(n) \leq \frac{4}{n}$ for $n \geq 2$. (ii) $p_0(n)=\frac{2+o(1)}{n}$.
\end{claim}
\noindent
{\bf Proof.} $X_{\pi}$ is disconnected iff there exists a $1 \leq k \leq n-1$
such that $1 \leq \pi(i) \leq n-k$ for all $k+1 \leq i \leq n$.
It follows that for $n\geq 5$
\begin{equation*}
\begin{split}
p_0(n) &\leq \frac{1}{n!}\sum_{k=1}^{n-1} k!(n-k)!
\leq \frac{2}{n}+\frac{4}{n(n-1)}+(n-5)\binom{n}{3}^{-1} \\
&\leq
\min\left\{\frac{4}{n}, \frac{2}{n}+ \frac{10}{n(n-1)}\right\}.
\end{split}
\end{equation*}
This proves the upper bounds in (i) and (ii). For the lower bounds note that if $n \geq 2$ then
\[
\pr [X_{\pi} {\rm ~is~disconnected~}] \geq \pr [\pi(1)=n {\rm ~or~} \pi(n)=1] =
\frac{2}{n}- \frac{1}{n(n-1)}.
\]
{\enp}

\noindent
{\bf Proof of Theorem \ref{t:maint}(i).} We argue by induction on $r$. The induction base $r=0$ follows from Claim \ref{c:pzero}(i).
Assume now that $r \geq 1$.
Applying (\ref{e:recursion}) and the induction hypothesis, it follows that
\begin{equation*}
\label{e:precu}
\begin{split}
&p_r(n) \leq 20 \left(1+\log n \right) \left( \frac{1}{n}+\sum_{k=0}^{n-3} \frac{p_{r-1}(n-2-k)}{k+2} \right) \\
&\leq
20 \log 3n \left(\frac{1}{n}+ 40^r \sum_{k=0}^{n-3}
 \frac{\log(3(n-k-2))^{2(r-1)}}{(k+2)(n-k-2)} \right) \\
&\leq  20 \log 3n \left(\frac{1}{n}+ 40^r (\log 3n)^{2(r-1)}\sum_{k=0}^{n-3}
\frac{1}{(k+2)(n-k-2)} \right) \\
&= 20 \log 3n \left(\frac{1}{n}+ \frac{40^r (\log 3n)^{2(r-1)}}{n} \sum_{k=0}^{n-3} \left(\frac{1}{k+2}+\frac{1}{n-k-2}\right) \right) \\
&= 20 \log 3n \left(\frac{1}{n}+ \frac{40^r (\log 3n)^{2(r-1)}}{n} \left( \sum_{k=2}^{n-1}\frac{1}{k}+\sum_{k=1}^{n-2} \frac{1}{k}\right) \right)\\
&\leq 20 \log 3n \left(\frac{1}{n}+ \frac{40^r (\log 3n)^{2(r-1)}}{n}(1+2 \log n) \right) \\
&\leq 40^{r+1} \frac{(\log 3n)^{2r}}{n}.
\end{split}
\end{equation*}
{\enp}

\section{The Lower Bound}
\label{s:lbnd}

The proof of Theorem \ref{t:maint}(ii) depends on a homotopy decomposition result for permutation complexes given in \cite{CLM20}.
For $\pi \in \sym_n$ and a subset $A\subset [n]$, let $\psi_A(\pi)\in \sym_{|A|}$ denote the permutation pattern of $\pi$
restricted to $A$, i.e. if $A=\{a_1,\ldots,a_k\}$ where $a_1<\cdots<a_k$, then $\psi_A(\pi)=\sigma^{-1}$
where $\sigma \in \sym_k$ satisfies $\pi(a_{\sigma(1)})<\cdots<\pi(a_{\sigma(k)})$.
For $1 \leq t  \leq n$ define
\[
\pi_t'=\psi_{[n]\setminus\{t\}}(\pi)\in \sym_{n-1}~~~,~~~
\pi_t''=\psi_{[n]\setminus[t]}(\pi) \in \sym_{n-t}.
\]
\noindent
When $t=n$, we view $\pi_t''$ as the empty permutation and define $X_{\pi_t''}=\emptyset$.
\ \\ \\
Let $\pi \in \sym_n$ and let $i=\pi^{-1}(1), j=\pi^{-1}(2)$.
As usual, let $\Sigma X=S^0*X$ denote the suspension of a complex $X$.
\begin{proposition}[\cite{CLM20}]
\label{p:homdec}
Let $\pi \in \sym_n$ and let $i=\pi^{-1}(1), j=\pi^{-1}(2)$. Then
\begin{itemize}
\item[(i)]
If $i<j$ then $X_{\pi} \simeq X_{\pi_{j}'}$.
\item[(ii)]
If $i>j$ then $X_{\pi} \simeq X_{\pi_i'} \vee \Sigma X_{\pi_{i}''}$.
\end{itemize}
\end{proposition}

\begin{proposition}
\label{p:lbound}
For $r \geq 1$ and $n \geq 2$
\begin{equation}
\label{e:lbound}
p_r(n) \geq \frac{1}{n}\left(1+\frac{1}{n-1}\sum_{i=1}^n (i-1) p_{r-1}(n-i)\right).
\end{equation}
\end{proposition}
\noindent
{\bf Proof.} Fix $1 \leq j< i \leq n$. By Proposition \ref{p:homdec}(ii)
\begin{equation}
\label{e:decij}
\begin{split}
&\pr\left[\pi \in \sym_n: \pi^{-1}(1)=i, \pi^{-1}(2)=j, X_{\pi} \text{~is~not~}r\text{-connected}\right] \\
&\geq \pr\left[\pi \in \sym_n: \pi^{-1}(1)=i, \pi^{-1}(2)=j, X_{\pi_{i}''}\text{~is~not~}(r-1)\text{-connected}\right] \\
&=\frac{p_{r-1}(n-i)}{n(n-1)}.
\end{split}
\end{equation}
Summing (\ref{e:decij}) over all $1 \leq j< i \leq n$ it follows that
\begin{equation*}
\label{e:feij}
\begin{split}
p_r(n) &\geq \sum_{1 \leq j < i \leq n} \pr\left[\pi \in \sym_n: \pi^{-1}(1)=i, \pi^{-1}(2)=j, X_{\pi} \text{~is~not~}r\text{-connected}\right] \\
&\geq \sum_{1 \leq j < i \leq n}\frac{p_{r-1}(n-i)}{n(n-1)}
=\frac{1}{n(n-1)}\sum_{1 \leq i \leq n}(i-1)p_{r-1}(n-i) \\
&=\frac{1}{n}\left(1+\frac{1}{n-1}\sum_{i=1}^{n-1} (i-1) p_{r-1}(n-i)\right).
\end{split}
\end{equation*}
{\enp}
\noindent
{\bf Proof of Theorem \ref{t:maint}(ii).} We argue by induction on $r$. The induction base $r=0$ follows from Claim \ref{c:pzero}(i).
Assume now that $r \geq 1$. Applying (\ref{e:lbound}) and the induction hypothesis we obtain
\begin{equation*}
\begin{split}
p_r(n) &\geq \frac{1}{n(n-1)}\sum_{i=1}^{n-1} (i-1) p_{r-1}(n-i) \\
&= \frac{1}{n(n-1)}\sum_{i=1}^{n-1} (n-1-i) p_{r-1}(i) \\
&\geq  \frac{1}{n(n-1)}\sum_{i=1}^{n-1} (n-1-i) \left(
\frac{1}{(r-1)!}\frac{(\log i)^{r-1}}{i}-c_{r-1}\frac{(\log i)^{r-2}}{i} \right) \\
&\geq  \frac{1}{n}\sum_{i=2}^{n-1} \left(\frac{(\log i)^{r-1}}{(r-1)!i}-c_{r-1}\frac{(\log i)^{r-2}}{i} \right)
-\frac{1}{(r-1)!n(n-1)}\sum_{i=2}^{n-1}(\log i)^{r-1} \\
&\geq \frac{1}{(r-1)! n} \int_1^n \frac{(\log t)^{r-1}dt}{t} - c_{r-1}\frac{(\log n)^{r-2}}{n} \sum_{i=2}^{n-1} \frac{1}{i} -\frac{1}{(r-1)!} \frac{(\log n)^{r-1}}{n} \\
&\geq \frac{1}{r!} \frac{(\log n)^r}{n} - \left(c_{r-1}+\frac{1}{(r-1)!}\right) \frac{(\log n)^{r-1}}{n} \\
&=\frac{1}{r!} \frac{(\log n)^r}{n} - c_r \frac{(\log n)^{r-1}}{n}.
\end{split}
\end{equation*}
{\enp}

\section{Concluding Remarks}
\label{s:con}

In this paper we studied the topological connectivity of the order complex $Y_q$ of a sequence of $n$ random points $q=(q_1,\ldots,q_n)$
in the unit square. More generally, let $Y_{d,q}$ denote the order complex determined by sequence of $n$ points
$q=(q_1,\ldots,q_n) \in \Omega_{d,n}=
\left([0,1]^d\right)^n$. Our work suggests a number of directions for further research.
\begin{itemize}
\item
One natural problem is to close the gap between the lower and upper bounds in Theorem \ref{t:maint}. It seems likely
that the lower bound $p_r(n)=\Omega\left(\frac{(\log n)^r}{n} \right)$ is closer to the truth.
\item
Theorem \ref{t:hom} asserts that $Y_{2,q}=Y_q$ is a wedge of spheres. Is there a simple
characterization of the homotopy type of $Y_{d,q}$ for general $d$?
\item
Let $p_{d,r}(n)$ denote the probability that $Y_{d,q}$ is $r$-connected for a random $q=(q_1,\ldots,q_n) \in \Omega_{d,n}$.
It would be interesting to extend Theorem \ref{t:maint} to general $d$.
Our proof of the upper bound $p_{2,r}(n)=p_r(n)=O\left(\frac{(\log n)^{2r}}{n}\right)$ does not use very specific properties of two dimensional posets (such as Theorem \ref{t:hom} and the homotopy decomposition in Proposition \ref{p:homdec}), and may thus be
relevant to higher dimensions as well.
\end{itemize}

\ \\ \\
{\bf Acknowledgement:} We would like to thank Russ Woodroofe for bringing reference \cite{PS18}
to our attention.

\end{document}